\documentclass[a4paper,11pt]{article}
\usepackage{graphicx}
\usepackage{psfrag}
\usepackage{epsfig}
\usepackage{subfigure}

\usepackage{amssymb}
\usepackage{psfrag}
\usepackage{amsmath,amsthm,amssymb,eucal}
\usepackage{mathrsfs}
\usepackage{bm,amsbsy}
\usepackage[all]{xy}

\usepackage{authblk}

\usepackage{hyperref}

\hypersetup{
    colorlinks,%
    citecolor=black,%
    filecolor=black,%
    linkcolor=black,%
    urlcolor=black
}

\numberwithin{equation}{section}

\newcommand{\newsec}[1]{\medskip {\bf \noindent #1. \quad \hspace{-0.6cm}}}

\newtheorem{thmm}{Theorem}[section]
\newtheorem{prop}[thmm]{Proposition}
\newtheorem{lemma}[thmm]{Lemma}

\begin{document}

\title{Reduction of Stokes-Dirac structures and gauge symmetry in port-Hamiltonian systems}

 \author[*]{Marko Seslija}
 \author[$\dag$]{Arjan van der Schaft}
 \author[*]{Jacquelien M.A. Scherpen}

 \affil[*]{Department of Discrete Technology and Production Automation, Faculty of Mathematics and Natural Sciences, University of Groningen, Nijenborgh 4, 9747 AG Groningen, The Netherlands, e-mail:~\{M.Seslija,\;J.M.A.Scherpen\}@rug.nl}
 \affil[$\dag$]{Johann Bernoulli Institute for Mathematics and Computer Science, University of Groningen, Nijenborgh 9, 9747 AG Groningen, The Netherlands, e-mail:~A.J.van.der.Schaft@rug.nl}

\date{May 5, 2012}

\maketitle

\begin{abstract}
Stokes-Dirac structures are infinite-dimensional Dirac structures defined in terms of differential forms on a smooth manifold with boundary. These Dirac structures lay down a geometric framework for the formulation of Hamiltonian systems with a nonzero boundary energy flow. Simplicial triangulation of the underlaying manifold leads to the so-called simplicial Dirac structures, discrete analogues of Stokes-Dirac structures, and thus provides a natural framework for deriving finite-dimensional port-Hamiltonian systems that emulate their infinite-dimensional counterparts. The port-Hamiltonian systems defined with respect to Stokes-Dirac and simplicial Dirac structures exhibit gauge and a discrete gauge symmetry, respectively. In this paper, employing Poisson reduction we offer a unified technique for the symmetry reduction of a generalized canonical infinite-dimensional Dirac structure to the Poisson structure associated with Stokes-Dirac structures and of a fine-dimensional Dirac structure to simplicial Dirac structures. We demonstrate this Poisson scheme on a physical example of the vibrating string.
\end{abstract}


\section{Introduction}
Geometric structures behind a variety of physical systems stemming from mechanics, electromagnetism and chemistry exhibit a remarkable unity enunciated by Dirac structures. The open dynamical systems defined with respect to these structures belong to the class of so-called port-Hamiltonian systems. These systems arise naturally from the energy-based modeling. Apart from offering a geometric content of Hamiltonian systems, Dirac structures supply a framework for modeling port-Hamiltonian systems as interconnected and constrained systems. From a network-modeling perspective, this means that port-Hamiltonian systems can be reticulated into a set of energy-storing elements, a set of energy-dissipating elements, and a set of energy port by which the interconnection of these blocks and environment is modeled. It is well-known that such a modeling strategy also utilizes control synthesis for these systems.

The port-Hamiltonian formalism transcends the lumped-parameter scenario and has been successfully applied to study of a number of distributed-parameter systems \cite{vdSM02, MacCDCP1}. The centrepiece of the efforts concerning infinite-dimensional case is the Stokes-Dirac structure. The canonical Stokes-Dirac structure is an infinite-dimensional Dirac structure defined in terms of differential forms on a smooth manifold with boundary. The Hamiltonian equations associated to this Dirac structure allow for non-zero energy exchange through the boundary.

Although the differential operator in the Stokes-Dirac structure, in the presence of nonzero boundary conditions, is not skew-symetric, it is possible to associate a (pseudo-)Poisson structure to the Stokes-Dirac structure \cite{vdSM02}. In the absence of algebraic constraints, the Stokes-Dirac structure specializes to a Poisson structure \cite{Dorfman}, and as such it can be derived through symmetry reduction from a canonical Dirac structure on the phase space \cite{JorisCDC}. How to conduct this reduction for the Poisson structure associated to the Stokes-Dirac structure on a manifold with boundary is the central theme of this paper.

\vspace{0.25cm}
{\textbf{Contribution and outline.}} This paper very closely follows \cite{JorisCDC}. The reduction scheme we are dealing with is the one from \cite{JorisCDC}, the only difference being in that we consider slightly augmented spaces in order to account for the behaviours associated with the boundary. The perspective as well as the notation in Section~2 are taken verbatim from \cite{JorisCDC}, but now for the generalized Dirac structures that allow for the formulation of open Hamiltonian systems. The proposed Poisson reduction is firstly applied in the reduction of a generalized cannonical Dirac structure to the Poisson structure associated with the Stokes-Dirac structure. In the context of dynamics, the canonical port-Hamiltonian systems are those defined as in \cite{Schlacher,Schoberl}, now only in the context of differential forms, while the reduced port-Hamiltonian systems are exactly those presented in \cite{vdSM02}. In the final section we demonstrate how this reduction applies to the Poisson reduction of the port-Hamiltonian systems on discrete manifolds \cite{SeslijaMathMod,SeslijaJGP}.

\section{Dirac structures and reduction}

Dirac structures were originally developed in \cite{Courant,Dorfman} as a generalization of symplectic and Poisson structures. The formalism of Dirac structure was employed as the geometric notion underpinning generalized power-conserving interconnections and thus allowing the Hamiltonian formulation of interconnected and constrained dynamical systems.

Let $Q$ be a manifold and define a pairing on $TQ \oplus T^\ast Q$ given by 
\[
	\left<\!\left<(v, \alpha), (w, \beta) \right>\!\right> = 
	\frac{1}{2} (\alpha(w) + \beta(v)).
\]
For a subspace $D$ of $TQ \oplus T^\ast Q$, we define the orthogonal complement $D^\perp$ as the space of all $(v, \alpha)$ such that $\left<\!\left<(v, \alpha), (w, \beta) \right>\!\right> = 0$ for all $(w, \beta)$.  A \emph{\textbf{Dirac structure}} is then a subbundle $D$ of $TQ \oplus T^\ast Q$ which satisfies $D = D^\perp$.

The notion of Dirac structures just entertained is suitable for the formulation of closed Hamiltonian systems, however, our aim is a treatment of open Hamiltonian systems in such a way that some of the external variables remain free port variables. 

Let $F$ be a linear vector space of external flows, with dual the space $F^\ast$ of external efforts. We deal with Dirac structures on the product space $Q\times F$. The pairing on $(TQ\times F )\oplus (T^\ast Q \times F)$ is given by
\begin{equation}\label{eq:bilformgen}
\begin{split}
	&\left<\!\left< \Big((v,f), ( \alpha,e)\Big), \left((w,\tilde{f}),  (\beta, \tilde{e})\right) \right>\!\right> \\
	&\quad = 
	\frac{1}{2} \left(\alpha(w) +e(\tilde{f})+ \beta(v)+\tilde{e}(f)\right).
	\end{split}
\end{equation}
A \textbf{\emph{generalized Dirac structure}} $D$ is a subbundle of $(TQ\times F )\oplus (T^\ast Q \times F)$ which is maximally isotropic under (\ref{eq:bilformgen}).

Canonical Dirac structure on $TQ\oplus T^\ast Q$ is considered to be a symplectic structure. However, in this paper we shall deal with slightly different canonical Dirac structures. To that end, let the map $\sharp: T^\ast Q \times F^\ast \rightarrow TQ \times F$ induces the \emph{\textbf{Poisson structure}} on $TQ \times F$.  The graph of $\sharp$ given by 
\begin{align}
	D_{T^\ast Q\times F^\ast} & := \{ (\sharp(\alpha, e), (\alpha,e)) :  \alpha \in T^\ast Q \,, e\in F^\ast \}
				\label{poisson}
\end{align}
is a Dirac structure. If the mapping $\sharp$ is symplectic on $TQ$, that is if $\sharp(\alpha,0)=0$ implies $\alpha=0$, the Dirac structure (\ref{poisson}) is the \textbf{\emph{generalized canonical Dirac structure}}.

There is a number of techniques for symmetry reduction of Dirac structures \cite{Guido, Yoshimura2}. The reduction considered in this paper is the Poisson reduction from \cite{JorisCDC}. For that purpose, let $G$ be a Lie group which acts on $Q$ from the right and assume that the quotient space $Q/G$ is again a manifold.  Denote the action of $g \in G$ on $q \in Q$ by $q \cdot g$ and the induced actions of $g \in G$ on $TQ\times F$ and $T^\ast Q\times F^\ast$ by $(v,f) \cdot g$ and $(\alpha,e) \cdot g$, for $v \in TQ$, $f\in F$, $\alpha \in T^\ast Q$, and $e\in F^\ast$.  The action on the $T^\ast Q \times F^\ast$ is defined by $\left< (\alpha, e) \cdot g, (v,f) \right> = \left< (\alpha,e), (v,f) \cdot g^{-1} \right>$.  In what follows, we will focus mostly on the reduced cotangent bundle $(T^\ast Q \times F^\ast)/G$.  In this paper, we will deal with the space denoted by $T^\ast Q/G\times F^\ast$.

Consider now the canonical Dirac structure on $T^\ast Q\times F^\ast$.  Let $\sharp : T^\ast Q \times F^\ast \rightarrow TQ \times F$ be the map \eqref{poisson} used in the definition of $D_{T^\ast Q\times F^\ast}$.  The reduced Dirac structure $D_{T^\ast Q/G\times F^\ast}$ on $T^\ast Q/G\times F^\ast$ can now be described as the graph of a reduced map $[\sharp] : T^\ast (T^\ast Q/G\times F^\ast) \rightarrow T(T^\ast Q/G\times F^\ast)$ defined as follows.

Let $\pi_G : T^\ast Q\times F^\ast \rightarrow T^\ast Q/G\times F^\ast$ be the quotient map and consider an element $(\rho, \pi, \rho_b)$ in $T^\ast Q\times F^\ast$.  The tangent map of $\pi_G$ at $(\rho, \pi, \rho_b)$ is denoted by $T_{(\rho, \pi, \rho_b)} \pi_G  :  T_{(\rho, \pi, \rho_b)} (T^\ast Q\times F^\ast) \rightarrow T_{(\rho, \pi, \rho_b)} (T^\ast Q/G\times F^\ast)$, and its dual by $T_{(\rho, \pi, \rho_b)}^\ast \pi_G : T^\ast_{\pi_G(\rho, \pi, \rho_b)} (T^\ast Q/G\times F^\ast) \rightarrow T^\ast_{\pi_G(\rho, \pi, \rho_b)} (T^\ast Q\times F^\ast)$.  The reduced map $[\sharp]$ now fits into the following extended commutative diagram

\begin{equation} \label{prescription}
\xymatrix{
T^\ast_{(\rho, \pi, \rho_b)}(T^\ast Q\times F^\ast) \ar[r]^\sharp & T_{(\rho, \pi, \rho_b)}(T^\ast Q\times F^\ast) 
	\ar[d]^{T_{(\rho, \pi, \rho_b)} \pi_G} \\
T^\ast_{\pi_G(\rho, \pi, \rho_b)}(T^\ast Q/G\times F^\ast) \ar[u]_{T^\ast_{(\rho, \pi, \rho_b)} \pi_G} \ar[r]_{[\sharp]} & T_{\pi_G(\rho, \pi, \rho_b)}(T^\ast Q/G\times F^\ast).
}\end{equation}

\medskip
\section{Constant Stokes-Dirac structures}

Throughout this paper, let $M$ be an oriented $n$-dimensional smooth manifold with a smooth $(n-1)$-dimensional boundary $\partial M$ endowed with the induced orientation, representing the space of spatial variables. By $\Omega^k(M)$, $k=0,1,\ldots,n$, denote the space of exterior $k$-forms on $M$, and by $\Omega^k(\partial M)$, $k=0,1,\ldots,n-1$, the space of $k$-forms on $\partial M$. A natural non-degenerative pairing between $\rho\in\Omega^k(M)$ and $\sigma\in\Omega^{n-k}(M)$ is given by $\langle \sigma|\rho \rangle=\int_M \sigma \wedge \rho$. Likewise, the pairing on the boundary $\partial M$ between $\rho \in \Omega^k(\partial M)$ and $\sigma \in \Omega^{n-k-1}(\partial M)$ is given by $\langle \sigma|\rho\rangle=\int_{\partial M} \sigma \wedge \rho$ \cite{vdSM02}.
 
\subsection{Stokes-Dirac structure}
For any pair $p,q$ of positive integers satisfying $p+q=n+1$, define the flow and effort linear spaces by
\begin{equation*}
\begin{split}
\mathcal{F}_{p,q}=\,&\Omega^p(M)\times \Omega^q(M)\times \Omega^{n-p}(\partial M)\,\\
\mathcal{E}_{p,q}=\,&\Omega^{n-p}(M)\times \Omega^{n-q}(M)\times \Omega^{n-q}(\partial M)\,.
\end{split}
\end{equation*}
The bilinear form on the product space $\mathcal{F}_{p,q}\times\mathcal{E}_{p,q}$ is
\begin{equation}\label{eq-aj9}
\begin{split}
\langle\!\langle (&\underbrace{f_p^1,f_q^1,f_b^1}_{\in \mathcal{F}_{p,q}},\underbrace{e_p^1,e_q^1,e_b^1}_{\in \mathcal{E}_{p,q}}), (f_p^2,f_q^2,f_b^2,e_p^2,e_q^2,e_b^2)\rangle\!\rangle\\
& =\int_M e_p^1\wedge f_p^2+e_q^1\wedge f_q^2+ e_p^2\wedge f_p^1+e_q^2\wedge f_q^1\\
&~~~+\int_{\partial M} e_b^1\wedge f_b^2+ e_b^2\wedge f_b^1\,.
  \end{split}
\end{equation}

\begin{thmm}[Stokes-Dirac structure \cite{vdSM02}]
Given linear spaces $\mathcal{F}_{p,q}$ and $\mathcal{E}_{p,q}$, and the bilinear form $\langle\!\langle, \rangle\!\rangle$, define the following linear subspace $\mathcal{D}$ of $\mathcal{F}_{p,q}\times \mathcal{E}_{p,q}$\vspace{-0.1cm}
\begin{equation}\label{eq-7Cont}
\begin{split}
\mathcal{D}=\big\{& (f_p,f_q,f_b,e_p,e_q,e_b)\in \mathcal{F}_{p,q}\times \mathcal{E}_{p,q}\big |\\
& \left(\begin{array}{c}f_p \\f_q\end{array}\right)=\left(\begin{array}{cc}0 & (-1)^{pq+1} {\mathrm{d}}\\ {\mathrm{d}} & 0\end{array}\right)\left(\begin{array}{c}e_p \\e_q\end{array}\right)\,,\\
& \left(\begin{array}{c}f_b \\e_b\end{array}\right)=\left(\begin{array}{cc}\mathrm{tr}&0\\0 & -(-1)^{n-q}\mathrm{tr}\end{array}\right)\left(\begin{array}{c}e_p \\e_q \end{array}\right)
\big \}\,,
 \end{split}
\end{equation}
where $\mathrm{d}$ is the exterior derivative and $\mathrm{tr}$ stands for a trace on the boundary $\partial M$. 
Then $\mathcal{D}=\mathcal{D}^\perp$, that is, $\mathcal{D}$ is a Dirac structure.
\end{thmm}

It is possible to associate a Poisson structure to the Stokes-Dirac structure $\mathcal{D}$. Here we just sketch the essence and refer the reader to \cite{vdSM02}.

The space of admissible efforts is $\mathcal{E}_{\mathrm{adm}}:=\{ e\in \mathcal{E}_{p,q}| \exists f \in \mathcal{F}_{p,q}~\textrm{such~that}~(f,e)\in \mathcal{D} \}$. The set of admissible mappings $\mathcal{K}_{\mathrm{adm}}:=\{ k: \mathcal{F}_{p,q} \rightarrow \mathbb{R} | \forall a\in \mathcal{F}_{p,q},\,  \exists e(k,a)\in \mathcal{E}_{\mathrm{adm}}~\textrm{such~that~for}~ \forall a\in \mathcal{F}_{p,q}~ k(a+\partial a)= k(a) +\langle e(k,a) | \partial a \rangle + O(\partial a) \}$. The set $\mathcal{K}_{\mathrm{adm}}$ consists of those functions $k: \Omega^p(M) \times \Omega^q(M) \times \Omega^{n-p}(\partial M) \rightarrow \mathbb{R}$ whose derivatives $\delta k(z) =( \delta_pk(z), \delta_q k(z), \delta_b k(z) ) \in \Omega^{n-p}(M) \times \Omega^{n-q}(M) \times \Omega^{n-q}(\partial M)$ satisfy $\delta_b k(z)= -(-1)^{n-q} \mathrm{tr}(\delta_q k(z))$. The Poisson bracket on $\mathcal{K}_{\mathrm{adm}}$ is given as
\begin{equation*}
\begin{split}
\!\{ k^1, k^2\}_\mathcal{D}=&\int_M ((\delta_p k^1) \wedge (-1)^r \mathrm{d}( (\delta_q k^2) + (\delta_q k^1) \wedge \mathrm{d} (\delta_p k^2) ) \\ &- \int_{\partial M} ( (-1)^{n-q} (\delta_q k^1 ) \wedge (\delta_p k^2) )\,.
\end{split}
\end{equation*}
Using Stokes' theorem, it follows that the bracket is skew-symmetric and that it satisfies the Jacobi identity: $\{ \{ k^1, k^2\}_\mathcal{D}, k^3\}_\mathcal{D}+ \{ \{ k^1
2, k^3\}_\mathcal{D}, k^1\}_\mathcal{D} + \{ \{ k^3, k^1\}_\mathcal{D}, k^2\}_\mathcal{D} =0$ for all $k^i \in \mathcal{K}_{\mathrm{adm}} $.

In this paper we will exclusively be dealing with Poisson and associated Poisson structures.

\subsection{Simplicial Dirac structures}

In the discrete setting, the smooth manifold $M$ is replaced by an $n$-dimensional well-centered oriented manifold-like simplicial complex $K$ \cite{SeslijaMathMod, SeslijaJGP}. The flow and the effort spaces will be the spaces of complementary primal and dual forms. The elements of these two spaces are paired via the discrete primal-dual wedge product. Let
\begin{equation*}
\begin{split}
\mathcal{F}_{p,q}^d&=\Omega_d^p(\star_\mathrm{i} K)\times \Omega_d^q( K)\times \Omega_d^{n-p}(\partial (K))\\
\mathcal{E}_{p,q}^d&=\Omega_d^{n-p}( K)\times \Omega_d^{n-q}(\star_\mathrm{i} K)\times \Omega_d^{n-q}(\partial (\star K))\,.
\end{split}
\end{equation*}

The primal-dual wedge product ensures a bijective relation between the primal and dual forms, between the flows and efforts. A natural discrete mirror of the bilinear form (\ref{eq-aj9}) is a symmetric pairing on the product space $\mathcal{F}_{p,q}^d\times \mathcal{E}_{p,q}^d$ defined by
\begin{equation}\label{eq:bild}
\begin{split}
\langle\!\langle (&\underbrace{\hat{f}_p^1,{f}_q^1,{f}_b^1}_{\in \mathcal{F}_{p,q}^d},\underbrace{{e}_p^1,\hat{e}_q^1,\hat{e}_b^1}_{\in \mathcal{E}_{p,q}^d}), (\hat{f}_p^2,{f}_q^2,{f}_b^2,{e}_p^2,\hat{e}_q^2,\hat{e}_b^2)\rangle\!\rangle_d\\
&= \langle {e}_p^1\wedge \hat{f}_p^2+\hat{e}_q^1\wedge {f}_q^2+ {e}_p^2\wedge \hat{f}_p^1+\hat{e}_q^2\wedge {f}_q^1,K \rangle\\
&~\;~+\langle \hat{e}_b^1\wedge {f}_b^2+ \hat{e}_b^2\wedge f_b^1,\partial K \rangle
  \,.
  \end{split}
\end{equation}
A discrete analogue of the Stokes-Dirac structure is the finite-dimensional Dirac structure constructed in the following theorem \cite{SeslijaMathMod}.

\vspace{0.1cm}
\begin{thmm}[Simplicial Dirac structure \cite{SeslijaMathMod}]\label{th:pdDirac}
Given linear spaces $\mathcal{F}_{p,q}^d$ and $\mathcal{E}_{p,q}^d$, and the bilinear form $\langle\!\langle, \rangle\!\rangle_d$. The linear subspace $\mathcal{D}_d\subset\mathcal{F}_{p,q}^d\times \mathcal{E}_{p,q}^d$ defined by
\begin{equation}\label{eq:Dir-prim-dual}
\begin{split}
& \mathcal{D}_{d}=\big\{ (\hat{f}_p, {f}_q, {f}_b,{e}_p,\hat{e}_q,\hat{e}_b)\in \mathcal{F}_{p,q}^d\times \mathcal{E}_{p,q}^d\big |\\
& \!\!\left(\begin{array}{c}\hat{f}_p \\ {f}_q \end{array} \right) =\left(\begin{array}{cc} 0  & (-1)^{r}  \mathbf{d}_\mathrm{i}^{n-q}\\  \mathbf{d}^{n-p} & 0\end{array} \right) \left(\begin{array}{c} {e}_p \\ \hat{e}_q\end{array} \right)+ (-1)^{r} \left(\begin{array}{c}  \mathbf{d}_\mathrm{b}^{n-q} \\ 0 \end{array}\right) \hat{e}_b\,,\\
& \begin{array}{c}~~~~~f_b \end{array}= ~(-1)^{p}\mathbf{tr}^{n-p}{e}_p \}\,,
 \end{split}
\end{equation}
with $r=pq+1$, is a Dirac structure with respect to the pairing $\langle\!\langle, \rangle\!\rangle_{d}$ .
\end{thmm}

The operators $\mathbf{d}^{n-p}$ is the discrete exterior operator mapping $\Omega_d^{n-p}( K)$ to $\Omega_d^{q}( K)$, and $\mathbf{d}_i^{n-q}$ is the dual discrete exterior derivative. Note that since $\mathbf{d}_\mathrm{i}^{n-q}=(-1)^q (\mathbf{d}^{n-p})^\textsc{t}$ and $\mathbf{d}_\mathrm{b}^{n-q}=(-1)^{n-p}(\mathbf{tr}^{n-p})^\textsc{t}$, the structure (\ref{eq:Dir-prim-dual}) is in fact a \emph{Poisson structure} on the state space $\Omega_d^p( \star_\mathrm{i}K)\times \Omega_d^q( K)$.

The simplicial Dirac structure (\ref{eq:Dir-prim-dual}) is used as \emph{terminus a quo} for the geometric formulation of spatially discrete port-Hamiltonian systems \cite{SeslijaJGP}.

\section{Reduction of Stokes-Dirac structure}
The configuration manifold is a vector space $Q := \Omega^k(M)$ with the tangent bundle $TQ = Q \times Q$ and the cotangent bundle $T^\ast Q = Q \times Q^\ast$, where $Q^\ast = \Omega^{n-k}(M)$. The space of the boundary flows $F$ will be an admissible subset of $\Omega^{n-k-1}(\partial M)$, while the space of the boundary efforts is $E:=F^\ast=\Omega^k(\partial M)$.

The tangent bundle $T ( T^\ast Q \times F^\ast)$ is isomorphic to $(Q \times Q^\ast \times F^\ast) \times (Q \times Q^\ast \times F^\ast)$, with a typical element denoted by $(\rho, \pi, \rho_b, \dot{\rho}, \dot{\pi}, \dot\rho_b)$, while $T^\ast (T^\ast Q\times F^\ast) = (Q \times Q^\ast \times F^\ast) \times (Q^\ast \times Q \times F)$, with a typical element denoted by $(\rho, \pi, \rho_b, e_\rho, e_\pi, e_b)$.  For the duality pairing between $T(T^\ast Q\times F^\ast)$ and $T^\ast (T^\ast Q\times F^\ast)$ we chose
\begin{equation} \label{origdual}
\begin{split}
	&\left<  (\rho, \pi, \rho_b, e_\rho, e_\pi, e_b), (\rho, \pi, \rho_b, \dot{\rho}, \dot{\pi}, \dot\rho_b) \right>
	\\&= \int_M ( e_\rho \wedge \dot{\rho} + e_\pi \wedge \dot{\pi}) +\int_{\partial M} (e_b \wedge \dot\rho_b + e_b \wedge \textrm{tr}\, \dot{\rho})\,.
	\end{split}
\end{equation}
The choice for this non-degenerate pairing will become clear later on.

\subsection{The symmetry group}
Let $G$ be an Abelian group of $(k-1)$-forms. For any $\alpha \in G$ and $\rho \in Q$, the group $G$ action on $Q$ is
\begin{equation} \label{addaction}
	\rho \cdot \alpha  = \rho + \mathrm{d} \alpha.
\end{equation}
This action of gauge group lifts to $TQ\times F$ and $T^\ast Q\times F^\ast$ as $
	 (\rho, \dot{\rho}, e_b) \cdot \alpha = (\rho + \mathrm{d} \alpha, \dot{\rho}, e_b)
	$ and $ 
	 (\rho, \pi, \rho_b) \cdot \alpha = (\rho + \mathrm{d} \alpha, \pi, \rho_b)
$
for $\alpha \in G$, $(\rho, \dot{\rho}, e_b) \in TQ\times F$ and $(\rho, \pi, \rho_b) \in T^\ast Q \times F^\ast$.

The elements of $Q/G$ are equivalence classes $[\rho]$ of $k$-forms up to exact forms, so that the exterior differential determines a well-defined map from $Q/G$ to ${\mathrm{d}} \Omega^k$, given by $[\rho] \mapsto {\mathrm{d}} \rho$.  If the $k$-th cohomology of $M$ vanishes, we have $Q / G = {\mathrm{d}} \Omega^k$.  Consequently, the quotient $(T^\ast Q/G \times F^\ast)$ is isomorphic to $Q/G \times Q^\ast \times F^\ast$, or explicitly
\[
	(T^\ast Q \times F^\ast)/G =  {\mathrm{d}} \Omega^k(M) \times \Omega^{n - k}(M) \times \Omega^{k}(\partial M). 
\]
The quotient map denoted as $\pi_G : T^\ast Q\times F^\ast \rightarrow (T^\ast Q)/G \times F^\ast$ is given by 
\begin{equation} \label{quotmap}
	\pi_G (\rho, \pi, \rho_b) = ({\mathrm{d}} \rho, \pi, \rho_b).
\end{equation}

Let a representative element of $T^\ast Q/G \times F^\ast$ be $(\bar{\rho}, \bar{\pi}, \bar{\rho}_b)$, with $\bar{\rho} \in {\mathrm{d}} \Omega^k(M)$, $\bar{\pi} \in \Omega^{n - k}(M)$ and $\bar \rho_b \in \Omega^{k}(\partial{M})$.  Elements of $T(T^\ast Q/G \times F^\ast)$ will be denoted by $(\bar \rho , \bar{\pi}, \bar \rho_b, \dot{\bar{\rho}}, \dot{\bar{\pi}}, \dot{\bar{\rho}}_b )$, while the elements of $T^\ast(T^\ast Q/G \times F^\ast)$ will be denoted by $(\bar{\rho}, \bar{\pi}, \bar{\rho}_b, \bar{e}_\rho, \bar{e}_\pi, \bar{e}_b)$.  For the duality pairing, we use
\begin{equation}
\begin{split}
	&\left<  (\bar{\rho}, \bar{\pi}, \bar{\rho_b}, \bar{e}_\rho, \bar{e}_\pi, \bar{e_b}), 
		(\bar{\rho}, \bar{\pi}, \bar{\rho_b}, \dot{\bar{\rho}}, \dot{\bar{\pi}}, \dot{\bar{\rho}}_b) \right>
	= \int_M ( \bar{e}_\rho \wedge \dot{\bar{\rho}} + 
		\bar{e}_\pi \wedge \dot{\bar{\pi}}  )+ \int_{\partial M} {\bar{e}_b} \wedge \dot{\bar{\rho}}_b. 
\end{split}\end{equation}
Whenever the base point $(\bar{\rho}, \bar{\pi}, \bar{\rho}_b)$ is clear from the context, we will denote $(\bar{\rho}, \bar{\pi}, \bar{\rho}_b, \dot{\bar{\rho}}, \dot{\bar{\pi}}, \dot{\bar{\rho}}_b)$ simply by $(\dot{\bar{\rho}}, \dot{\bar{\pi}}, \dot{\bar{\rho}}_b)$, and similarly for $(\bar{\rho}, \bar{\pi}, \bar{\rho}_b, \bar{e}_\rho, \bar{e}_\pi, \bar{e}_b)$.

\subsection{The reduced Dirac structure}
The generalized canonical Dirac structure is a Poisson structure induced by the linear maping $\sharp : T^\ast( T^\ast Q \times F^\ast) \rightarrow T(T^\ast Q \times F^\ast)$ given by 
\begin{equation} \label{canonsymp}
	\sharp ( \rho, \pi, \rho_b, e_\rho, e_\pi, e_b) = (\rho, \pi, \rho_b, e_\pi, -(-1)^{k(n-k)}e_\rho, -\textrm{tr\,} e_\pi).
\end{equation}

In order to obtain the reduced Poisson structure from the canonical Dirac structure (\ref{canonsymp}), we need to specify what are the operators $T \pi_G$ and $T^\ast \pi_G$ in the diagram (\ref{prescription}). The space $F$ is the set of admissible forms $\Omega^{n-k-1}(\partial M)$ that are the traces of $(\mathrm{d}\Omega^{k})^\ast$, as will be made clear in Lemma~\ref{lemmaTTs}. Consider an element $(\rho, \pi, \rho_b) \in T^\ast Q\times F^\ast$, and we recall that $\pi_G(\rho, \pi, \rho_b) = ({\mathrm{d}} \rho, \pi, \rho_b)$.  Let $T_{(\rho, \pi, \rho_b)} \pi_G : T_{(\rho, \pi, \rho_b)} (T^\ast Q \times F^\ast) \rightarrow T_{({\mathrm{d}}\rho, \pi, \rho_b)} (T^\ast Q/G \times F^\ast)$ be the tangent map to $\pi_G$ at $(\rho, \pi, \rho_b)$ and consider the adjoint map $T^\ast_{(\rho, \pi, \rho_b)} \pi_G :  T^\ast_{({\mathrm{d}}\rho, \pi, \rho_b)} (T^\ast Q/G \times F^\ast) \rightarrow  T^\ast_{(\rho, \pi, \rho_b)} (T^\ast Q \times F^\ast)$.

\begin{lemma}\label{lemmaTTs}
The tangent and cotangent maps $T_{(\rho, \pi, \rho_b)} \pi_G$ and $T^\ast_{(\rho, \pi, \rho_b)} \pi_G$ are given by 
	\begin{equation} \label{tangent}
		T_{(\rho, \pi, \rho_b)} \pi_G(\rho, \pi, \rho_b, \dot{\rho}, \dot{\pi}, \dot{\rho}_b)
			= ({\mathrm{d}} \rho, \pi, \rho_b, {\mathrm{d}}\dot{\rho}, \dot{\pi}, \dot{\rho}_b)
	\end{equation}
and 
	\begin{equation} \label{cotangent}
	\begin{split}
		&T^\ast_{(\rho, \pi, \rho_b)} \pi_G({\mathrm{d}} \rho, \pi, \rho_b, \bar{e}_\rho, \bar{e}_\pi, -(-1)^{n-k}{\mathrm{tr}}\,\bar e_\rho) \\ &\quad=
			(\rho, \pi, \rho_b, (-1)^{n-k} {\mathrm{d}} \bar{e}_\rho, \bar{e}_\pi, -(-1)^{n-k}{\mathrm{tr}}\, \bar{e}_\rho).
	\end{split}\end{equation}
\end{lemma}

\medskip
\begin{proof}
The expression \eqref{tangent} for $T_{(\rho, \pi, \rho_b)} \pi_G$ follows from (\ref{quotmap}).  To prove \eqref{cotangent}, we let 
$(\dot{\rho}, \dot{\pi}, \dot{\rho}_b) \in T_{(\rho, \pi, \rho_b)} (T^\ast Q \times F^\ast)$ and consider 
\begin{align*}
\hspace{-0.2cm}&\left< T^\ast_{(\rho, \pi, \rho_b)} \pi_G(\bar{e}_\rho, \bar{e}_\pi, -(-1)^{n-k}\bar{e}_\rho), 
	(\dot{\rho}, \dot{\pi}, \dot{\rho}_b) \right> \\ & = \left< (\bar{e}_\rho, \bar{e}_\pi, -(-1)^{n-k} \mathrm{tr}\,\bar{e}_\rho), 
	T_{(\rho, \pi, \rho_b)} \pi_G(\dot{\rho}, \dot{\pi}, \dot{\rho}_b) \right> \\
	& = \left< (\bar{e}_\rho, \bar{e}_\pi, -(-1)^{n-k} \mathrm{tr}\, \bar{e}_\rho), 
	({\mathrm{d}}\dot{\rho}, \dot{\pi}, \dot{\rho}_b) \right>.
\end{align*}
Applying Stokes' theorem, we have
\begin{equation*}
\begin{split}
&\left< (\bar{e}_\rho, \bar{e}_\pi, -(-1)^{n-k} \mathrm{tr}\, \bar{e}_\rho), 
	({\mathrm{d}}\dot{\rho}, \dot{\pi}, \dot{\rho}_b) \right>\\
& = \int_M ( \bar{e}_\rho \wedge {\mathrm{d}}\dot{\rho} + \bar{e}_\pi \wedge\dot{\pi}) \\
&\quad+\int_M (-(-1)^{n-k} \mathrm{tr}\, \bar{e}_\rho \wedge \rho_b -(-1)^{n-k} \mathrm{tr}\, \bar{e}_\rho \wedge \mathrm{tr}\, \dot{\rho}) \\
& = \int_M ( (-1)^{n-k}{\mathrm{d}}\bar{e}_\rho \wedge \dot{\rho} + \bar{e}_\pi \wedge\dot{\pi}) \!- \!\! \int_{\partial M} (-1)^{n-k} \mathrm{tr}\, \bar{e}_\rho \wedge \mathrm{tr}\,\dot{\rho}\,.
\end{split}
\end{equation*}
Thus, $T^\ast_{(\rho, \pi, \rho_b)} \pi_G(\bar{e}_\rho, \bar{e}_\pi, -(-1)^{n-k} \mathrm{tr}\, \bar{e}_\rho) = ( (-1)^{n-k}{\mathrm{d}}\bar{e}_\rho,  \bar{e}_\pi, -(-1)^{n-k} \mathrm{tr}\, \bar{e}_\rho)$.
\end{proof}
\medskip

As in the case of a boundaryless manifold \cite{JorisCDC}, the reduced Poisson structure in \eqref{prescription} is given by
\[
	[\sharp]_{({\mathrm{d}} \rho, \pi, \rho_b)} = 
		T_{(\rho, \pi, \rho_b)} \pi_G \circ \sharp \circ T^\ast_{({\mathrm{d}} \rho, \pi, \rho_b)} \pi_G
\]
for all $({\mathrm{d}} \rho, \pi, \rho_b) \in T^\ast Q/G \times F^\ast$. 

\begin{thmm}The reduced Poisson structure is given by
\begin{equation} \label{redpoisson}
	[\sharp](\bar{e}_\rho, \bar{e}_\pi, -(-1)^{n-k} \mathrm{tr}\, \bar{e}_\rho) = 
		( {\mathrm{d}} \bar{e}_\pi, -(-1)^{n(k+1)} {\mathrm{d}} \bar{e}_\rho, -\mathrm{tr}\, \bar{e}_\pi).
\end{equation}
\end{thmm}

\newsec{Relation to the Stokes-Dirac structure}
The matrix form of the reduced Poisson structure is
\begin{equation} \label{matform1}
\hspace{-0.3cm}	\begin{pmatrix}
		\dot{\bar{\rho}} \\
		\dot{\bar{\pi}} \\
		\dot{\bar{\rho}}_b
	\end{pmatrix}
	= 
	\begin{pmatrix}
		0 & {\mathrm{d}} & 0\\
		-(-1)^{n(k+1)} {\mathrm{d}} & 0 & 0\\
		0 & -\mathrm{tr}\;\; & 0  
	\end{pmatrix}
	\begin{pmatrix}
		\bar{e}_\rho \\
		\bar{e}_\pi \\
		(-1)^{n-k}\mathrm{tr}\,\bar{e}_\rho
	\end{pmatrix}.
\end{equation}
The sign convention in (\ref{matform1}) and \cite{vdSM02} is not the same. To match the signs we introduce new \emph{flow variables} $f_p, f_q, f_b$ and \emph{effort variables} $e_p, e_q, e_b$ defined as $ e_p = \bar{e}_\rho,  e_q = (-1)^r \bar{e}_\pi, 	f_p = \dot{\bar{\rho}},  f_q = (-1)^{n(k+1)+1} \dot{\bar{\pi}},  f_b= -(-1)^r \dot{\bar{\rho}}_b $, where $p = k+1 $, $q = n-k$, and $r = pq + 1$.  With this choice of signs, \eqref{matform1} becomes 
\begin{equation}
\begin{split}
	\begin{pmatrix}
		f_p \\
		f_q
	\end{pmatrix}
	&= 
	\begin{pmatrix}
		0 & (-1)^r{\mathrm{d}}  \\
		{\mathrm{d}} & 0   
	\end{pmatrix}
	\begin{pmatrix}
		e_p \\
		e_q 
	\end{pmatrix}\\
		f_b&= \mathrm{tr}\,e_q\\
	(-1)^{n-k} \mathrm{tr}\,\bar{e}_\rho&=e_b\,.
\end{split}
\end{equation}

Here, it is important to point out that the boundary effort $e_b$, unlike in the case of the Stokes-Dirac structure, does not follow from the associate Poisson structure, but rather belongs to the set of admissible derivatives of the flow restricted to the boundary.

\medskip
\section{Symmetry in port-Hamiltonian systems}
Let $t \mapsto (\alpha_{\dot{\rho}}, \alpha_{\dot{\pi}}) \in \Omega^k(M) \times \Omega^{n-k}(M)$ be a time function, and let the Hamiltonian be
\[
H(\alpha_{{\dot{\rho}}}, \alpha_{\dot{\pi}})=\int_M \mathcal{H}(\mathrm{d}\alpha_{{\dot{\rho}}}, \alpha_{\dot{\pi}})\,.
\]
It follows that at any time instance $t\in \mathbb{R}$
\[
\frac{\mathrm{d} H}{\mathrm{d} t}=\int_M \frac{\delta H}{\delta \alpha_{{\dot{\rho}}}} \wedge \frac{\partial \alpha_{{\dot{\rho}}}}{\partial t} + \frac{\delta H}{\delta \alpha_{{\dot{\pi}}}} \wedge \frac{\partial \alpha_{{\dot{\pi}}}}{\partial t} +\int_{\partial M} \frac{\partial \mathcal{H}}{\partial (\mathrm{d} \alpha_{{\dot{\rho}}})} \wedge \frac{\partial \alpha_{{\dot{\rho}}}}{\partial t} \,.
\]
The differential forms $\frac{\partial \alpha_{{\dot{\rho}}}}{\partial t}, \frac{\partial \alpha_{{\dot{\pi}}}}{\partial t}$ represent the generalized velocities of the energy variables $\alpha_{{\dot{\rho}}}, \alpha_{{\dot{\pi}}}$. The connection with the canonical Dirac structure is made by setting the flows
\[
\dot{\rho}= -\frac{\partial \alpha_{\dot{\rho}}}{\partial t}\, \quad \dot{\pi}= -\frac{\partial \alpha_{\dot{\pi}}}{\partial t}\,,
\]
and the efforts
\[
e_\rho = \frac{\delta H}{\delta \alpha_{{\dot{\rho}}}}\,, \quad e_\pi = \frac{\delta H}{\delta \alpha_{{\dot{\pi}}}}\,.
\]

The \emph{\textbf{canonical distributed-parameter port-Hamilto-nian system}} on an $n$-dim-ensional manifold, with the state space $\Omega^k(M) \times \Omega^{n-k}(M)$, the Hamiltonian $H$ and the canonical Dirac structure (\ref{canonsymp}), is given as
\begin{equation}\label{canpHsys}
\begin{split}
\begin{pmatrix}
		- \frac{\partial \alpha_{\dot{\rho}}}{\partial t} \\
		- \frac{\partial \alpha_{\dot{\pi}}}{\partial t} 
	\end{pmatrix}
&=
\begin{pmatrix}
		0 & 1 \\
		-(-1)^{k(n-k)} & 0
	\end{pmatrix}
\begin{pmatrix}
		\frac{\delta H}{\delta  \alpha_{\dot{\rho}} } \\
		\frac{\delta H}{\delta  \alpha_{\dot{\pi}}}
\end{pmatrix}\\
\begin{pmatrix}
		f_b \\
		e_b 
	\end{pmatrix}
&=
\begin{pmatrix}
		0 & -\mathrm{tr} \\
		-\mathrm{tr} & 0
	\end{pmatrix}
\begin{pmatrix}
		\frac{\partial  \mathcal{H}}{\partial (\mathrm{d}  \alpha_{\dot{\rho}} )} \\
		\frac{\delta H}{\delta  \alpha_{\dot{\pi}}}
\end{pmatrix}\,.
\end{split}
\end{equation}

\begin{prop}
For the port-Hamiltonian system (\ref{canpHsys}) the following property
\[
\frac{\mathrm{d} H}{\mathrm{d} t}=\int_{\partial M} e_b \wedge f_b\,
\]
expresses the fact that the increase in energy on the domain $M$ is equal to the power supplied to the system through the boundary $\partial M$.
\end{prop}

\subsection{The reduced port-Hamiltonian systems}
The Hamiltonian $H$ is invariant if a spatially independent $k$-form is added to $\alpha_{\dot{\rho}}$, thus the Poisson reduction is applicable. Let the reduced field be $\bar\alpha_{\dot{\rho}}:= \mathrm{d} \alpha_{\dot{\rho}}$, then the reduced Hamiltonian is
\[
H(\bar\alpha_{{\dot{\rho}}}, \alpha_{\dot{\pi}})=\int_M \mathcal{H}(\bar\alpha_{{\dot{\rho}}}, \alpha_{\dot{\pi}})\,.
\]
The port-Hamiltonian system with respect to the reduced Poisson structure is
\begin{equation}\label{redpHsys}
\begin{split}
\begin{pmatrix}
		- \frac{\partial \bar\alpha_{\dot{\rho}}}{\partial t} \\
		- \frac{\partial \alpha_{\dot{\pi}}}{\partial t} 
	\end{pmatrix}
&=
\begin{pmatrix}
		0 & \mathrm{d} \\
		-(-1)^{n(k+1)} \mathrm{d}& 0
	\end{pmatrix}
\begin{pmatrix}
		\frac{\delta H}{\delta  \bar{\alpha}_{\dot{\rho}} } \\
		\frac{\delta H}{\delta  \alpha_{\dot{\pi}}}
\end{pmatrix}\\
\begin{pmatrix}
		f_b \\
		e_b 
	\end{pmatrix}
&=
\begin{pmatrix}
		0 & -\mathrm{tr} \\
		-(-1)^{n-k}\mathrm{tr} & 0
	\end{pmatrix}
\begin{pmatrix}
		\frac{\delta H}{\delta  \bar{\alpha}_{\dot{\rho}} } \\
		\frac{\delta H}{\delta  \alpha_{\dot{\pi}}}
\end{pmatrix}\,.
\end{split}
\end{equation}
This is precisely the port-Hamiltonian system given in \cite{vdSM02}.

We will show how the general considerations of the reduction of port-Hamiltonian systems apply to a physical example of the vibrating string.

\subsection{Vibrating string}
Consider an elastic string of length $l$, elasticity modulus $T$, and mass density $\mu$, subject to traction forces at its ends. The underlying manifold is the segment $M=[0,l]\subset \mathbb{R}$, with coordinate $z$.

Under the assumption of linear elasticity, the Hamiltonian is given by
\[
H(u,p)=\int_M \mathcal{H}(u,p)= \frac{1}{2} \int_M (\mu^{-1} p \wedge * p + T \mathrm{d} u \wedge * \mathrm{d }u)\,,
\]
where $p\in \Omega^1(M)$ is the momentum conjugate to the displacement $u\in \Omega^0(M)$, and $*$ is the Hodge star.

The canonical Hamiltonian equations are
\begin{equation}\label{canonvib}
\begin{split}
\begin{pmatrix}
		\frac{\partial u}{\partial t} \\
		\frac{\partial p}{\partial t}
	\end{pmatrix}
&=
\begin{pmatrix}
		0 & 1 \\
		-1 & 0
	\end{pmatrix}
\begin{pmatrix}
		\frac{\delta H}{\delta u} \\
		\frac{\delta H}{\delta p}
\end{pmatrix}\\
f_b &= \mathrm{tr} (* \mu^{-1} p)\\
e_b &= \mathrm{tr} \left(\frac{\partial \mathcal{H}}{\partial (\mathrm{d} u)}\right)\,,
\end{split}
\end{equation}
or component-wise
\begin{equation*}
\begin{split}
\frac{\partial u}{\partial t} &= * \mu^{-1} p \\
\frac{\partial p}{\partial t} &= \mathrm{d} (* T \, \mathrm{d} u)\\
f_b &= \mathrm{tr} (* \mu^{-1} p)\\
e_b &= \mathrm{tr} (* T \,\mathrm{d} u)\,.
\end{split}
\end{equation*}
The Hamiltonian formulation (\ref{canonvib}) is identical to the formulation of the heavy chain system in \cite{Schoberl}.

The energy balance for the vibrating string is
\begin{equation*}
\begin{split}
\frac{\mathrm{d}H}{\mathrm{d} t}&=\int_M \frac{\delta H}{\delta u} \wedge \frac{\partial u}{\partial t} + \frac{\delta H}{\delta p} \wedge \frac{\partial p}{\partial t} + \int_{\partial M} \frac{\partial \mathcal{H}}{\partial (\mathrm{d} u)} \wedge \frac{\partial u}{\partial t}\\
&=\int_M -\mathrm{d} (* T \, \mathrm{d} u)\wedge *\mu^{-1} p + *\mu^{-1} p \wedge \mathrm{d} (* T \, \mathrm{d} u) \\
&\quad+ \int_{\partial M} * \mu^{-1} p \wedge * T \,\mathrm{d} u\\
& = \int_{\partial M} * \mu^{-1} p \wedge * T \,\mathrm{d} u = \int_{\partial M} e_b \wedge f_b\,.
\end{split}
\end{equation*}

The Hamiltonian is invariant if a time function is added to $u$. The potential energy can be expressed in terms of the strain $\alpha= \mathrm{d} u$ so that the reduced Hamiltonian is given by
\[
H_r(\alpha,p)=\int_M \mathcal{H}_r(u,p)= \frac{1}{2} \int_M (\mu^{-1} p \wedge * p + T \alpha \wedge * \alpha)\,.
\]
The Hamiltonian equations of the vibrating string now read as
\begin{equation}
\begin{split}
\begin{pmatrix}
		\frac{\partial \alpha}{\partial t} \\
		\frac{\partial p}{\partial t}
	\end{pmatrix}
&=
\begin{pmatrix}
		0 & \mathrm{d} \\
		\mathrm{d} & 0
	\end{pmatrix}
\begin{pmatrix}
		\frac{\delta H_r}{\delta \alpha} \\
		\frac{\delta H_r}{\delta p}
\end{pmatrix}\\
\begin{pmatrix}
		f_b \\
		e_b
	\end{pmatrix}
&=
\begin{pmatrix}
		0 & \mathrm{tr} \\
		\mathrm{tr} & 0
	\end{pmatrix}
\begin{pmatrix}
		\frac{\delta H_r}{\delta \alpha} \\
		\frac{\delta H_r}{\delta p}
\end{pmatrix}\,.
\end{split}
\end{equation}
These are the equations that correspond to the formulation of the vibration string system with respect to the Stokes-Dirac structure as is given in \cite{vdSM02}.

\medskip
\section{Symmetry reduction in discrete setting}
In the discrete world, the configuration space is the set of primal discrete forms $Q=\Omega^k(K)$ with the dual $Q^\ast=\Omega^{n-k}(\star_i K)$.  The space of the boundary efforts is $E=F^\ast=\Omega^{k}(\partial(K))$, and the space of the boundary flows is $F=\Omega^{n-k-1}(\partial (\star K))$. 

For the duality pairing between $T(T^\ast Q\times F^\ast)$ and $T^\ast (T^\ast Q\times F^\ast)$ we choose
\begin{equation} \label{origdual}
\begin{split}
	\hspace{-0.5cm}&\left<  (\rho, \pi, \rho_b, e_\rho, e_\pi, e_b), (\rho, \pi, \rho_b, \dot{\rho}, \dot{\pi}, \dot{\rho_b}) \right>
	= \int_M ( e_\rho \wedge \dot{\rho} + e_\pi \wedge \dot{\pi}) +\int_{\partial M} e_b \wedge \dot{\rho}_b\,,
\end{split}
\end{equation}
where $\wedge$ is the primal-dual wedge product.

The generalized canonical Dirac structure is a Poisson structure induced by the linear mapping $\sharp : T^\ast( T^\ast Q \times F^\ast) \rightarrow T(T^\ast Q \times F^\ast)$ given by 
\begin{equation} \label{canonsympd}
\begin{split}
	\!\!\!\!\!\!\!\!\! &\sharp  ( \rho, \pi, \rho_b, e_\rho, e_\pi, e_b) = (\rho, \pi, \rho_b, e_\pi, -(-1)^{k(n-k)}(e_\rho +\mathbf{d}_b^{n-k-1} \bar{e}_b), -{\mathbf{tr}}^k e_\pi).
	\end{split}
\end{equation}

The group $G$ that acts on $Q$ is described by the following action
\[
\alpha \cdot \rho = \rho + \mathbf{d}^{k-1} \alpha 
\]
for $\alpha\in G$ and $\rho\in Q$, where $ \mathbf{d}^{k-1} $ is the discrete exterior derivative.

The quotient is $(T^\ast Q/G \times F^\ast)= \mathbf{d}^{k}\Omega^k(K)\times \Omega^{n-k}(\star_i K) \times \Omega^{k}(\partial(K))$.

As in the continuous setting, the quotient map denoted as $\pi_G : T^\ast Q\times F^\ast \rightarrow (T^\ast Q)/G \times F^\ast$ is given by 
\begin{equation} \label{quotmap}
	\pi_G (\rho, \pi, \rho_b) = (\mathbf{d}^k \rho, \pi, \rho_b).
\end{equation}

For the duality pairing between $T^\ast(T^\ast Q/G \times F^\ast)$ and $T(T^\ast Q/G \times F^\ast)$, we take
\begin{equation*}
\begin{split}
	&\left<  (\bar{\rho}, \bar{\pi}, \bar{\rho_b}, \bar{e}_\rho, \bar{e}_\pi, \bar{e_b}), 
		(\bar{\rho}, \bar{\pi}, \bar{\rho_b}, \dot{\bar{\rho}}, \dot{\bar{\pi}}, \dot{\bar{\rho_b}}) \right>
	= \int_M ( \bar{e}_\rho \wedge \dot{\bar{\rho}} + 
		\bar{e}_\pi \wedge \dot{\bar{\pi}}  )+ \int_{\partial M} {\bar e_b} \wedge \dot{\bar \rho}_b. 
\end{split}
\end{equation*}
As before, whenever the base point $(\bar{\rho}, \bar{\pi}, \bar{\rho_b})$ is clear, we will denote $(\bar{\rho}, \bar{\pi}, \bar{\rho_b}, \dot{\bar{\rho}}, \dot{\bar{\pi}}, \dot{\bar \rho}_b)$ simply by $(\dot{\bar{\rho}}, \dot{\bar{\pi}}, \dot{\bar{\rho_b}})$, and similarly for $(\bar{\rho}, \bar{\pi}, \bar{\rho_b}, \bar{e}_\rho, \bar{e}_\pi, \bar{e}_b)$.

\begin{lemma}  
The tangent and cotangent maps $T_{(\rho, \pi, \rho_b)} \pi_G$ and $T^\ast_{(\rho, \pi, \rho_b)} \pi_G$ are given by 
	\begin{equation} \label{tangentd}
		T_{(\rho, \pi, \rho_b)} \pi_G(\rho, \pi, \rho_b, \dot{\rho}, \dot{\pi}, \dot{\rho}_b)
			= (\mathbf{d}^k \rho, \pi, \rho_b, \mathbf{d}^k\dot{\rho}, \dot{\pi}, \dot{\rho}_b)
	\end{equation}
and 
	\begin{equation} \label{cotangentd}
	\begin{split}
		&T^\ast_{(\rho, \pi, \rho_b)} \pi_G(\mathbf{d}^k \rho, \pi, \rho_b, \bar{e}_\rho, \bar{e}_\pi, \bar{e}_b) =
			(\rho, \pi, \rho_b, (-1)^{n-k} \mathbf{d}_i^{n-k-1} \bar{e}_\rho, \bar{e}_\pi, \bar{e}_b).
	\end{split}\end{equation}
\end{lemma}

\begin{thmm}[Reduced simplicial Dirac structure]The reduced simplicial Poisson structure is given by
$
	[\sharp](\bar{e}_\rho, \bar{e}_\pi, \bar{e}_b) = 
		( \mathbf{d}^k \bar{e}_\pi, -(-1)^{n(k+1)} ((-1)^{n-k} \mathbf{d}_i^{n-k-1} \bar{e}_\rho + \mathbf{d}_b^{n-k-1}\bar{e}_b), -{\mathbf{tr}}^k \bar{e}_\pi)
$.
\end{thmm}

\subsection*{Port-Hamiltonian systems on a simplicial complex}
The canonical port-Hamiltonian system with respect to the canonical Dirac structure is
\begin{equation}
\begin{split}
-\frac{\partial \alpha_{\dot\rho}}{\partial t} &= \frac{\partial H}{\partial \alpha_{\dot\pi}}(\alpha_{\dot\rho}, \alpha_{\dot\pi})\\
-\frac{\partial \alpha_{\dot\pi}}{\partial t} &= (-1)^{k(n-k)}\left( \frac{\partial H}{\partial \alpha_{\dot\rho}}(\alpha_{\dot\rho}, \alpha_{\dot\pi}) + \mathbf{d}_b^{n-k-1} \bar{e}_b \right)\\
\dot{\rho}_b &=-\mathbf{tr}^k \frac{\partial H}{\partial \alpha_{\dot\pi}}(\alpha_{\dot\rho}, \alpha_{\dot\pi})
\end{split}
\end{equation}
The rank of the underlying Poison structure is the rank of the symplectic phase space $\Omega^k(K) \times \Omega^{n-k}(\star_i K)$.

The canonical Hamiltonian $(\alpha_{\dot\rho}, \alpha_{\dot\pi}) \mapsto H(\alpha_{\dot\rho}, \alpha_{\dot\pi})$ can be expressed as
\begin{equation}
H(\alpha_{\dot\rho}, \alpha_{\dot\pi}):= H_r(\mathbf{d}^k \bar\alpha_{\dot{\rho}}, \alpha_{\dot\pi})\,.
\end{equation}

The reduced port-Hamiltonian equations assume the following form
\begin{equation*}
\begin{split}
-\frac{\partial \dot\alpha_{\dot{\rho}}}{\partial t} &=-\mathbf{d}^k\frac{\partial \alpha_{\dot{\rho}}}{\partial t}= \mathbf{d}^k\frac{\partial H}{\partial \alpha_{\dot\pi}}(\alpha_{\dot\rho}, \alpha_{\dot\pi})=\mathbf{d}^k \frac{\partial H_r}{\partial \alpha_{\dot\pi}}(\bar\alpha_{\dot{\rho}}, \alpha_{\dot\pi})\\
-\frac{\partial \alpha_{\dot\pi}}{\partial t} &= (-1)^{k(n-k)}\left( \frac{\partial H}{\partial \alpha_{\dot\rho}}(\alpha_{\dot\rho}, \alpha_{\dot\pi}) + \mathbf{d}_b^{n-k-1} \bar{e}_b \right) \\
&= (-1)^{k(n-k)}\left( (-1)^{n-k} \mathbf{d}_i^{n-k-1} \frac{\partial H_r}{\partial \bar \alpha_{\dot{\rho}}}(\bar \alpha_{\dot{\rho}}, \alpha_{\dot\pi}) + \mathbf{d}_b^{n-k-1} \bar{e}_b \right)\\
\dot{\rho}_b &=-\mathbf{tr}^k \frac{\partial H}{\partial \alpha_{\dot\pi}}(\alpha_{\dot\rho}, \alpha_{\dot\pi})=-\mathbf{tr}^k \frac{\partial H_r}{\partial \alpha_{\dot\pi}}(\bar\alpha_{\dot{\rho}}, \alpha_{\dot\pi})\,.
\end{split}
\end{equation*}
This is precisely the port-Hamiltonian system on a simplicial manifold as presented in \cite{SeslijaMathMod, SeslijaJGP}.

\section{Final remark}
This paper addresses the issue of the symmetry reduction of the generalized canonical Dirac structure to the Poisson structure associated with the Stokes-Dirac structure. The open avenue for the future work is to find a reduction procedure that would directly lead to the Stokes-Dirac structure. 

\section*{Acknowledgments}
The first author expresses his warmest thanks to Joris Vankerschaver for his hospitality during the author's stay in California.

\end{document}